\documentclass[10pt,oneside,reqno]{amsart}
\usepackage{txfonts}
\usepackage{bbm}
\usepackage{amsmath}
\usepackage{graphicx}
\usepackage{mathrsfs}
\usepackage{stmaryrd}
\usepackage{amsfonts}
\usepackage{enumerate,amsmath,amssymb,amsthm}
\newcommand{\be}{\begin{eqnarray}}
\newcommand{\ee}{\end{eqnarray}}
\newcommand{\ce}{\begin{eqnarray*}}
\newcommand{\de}{\end{eqnarray*}}
\newtheorem{theorem}{Theorem}[section]
\newtheorem{lemma}[theorem]{Lemma}
\newtheorem{remark}[theorem]{Remark}
\newtheorem{definition}[theorem]{Definition}
\newtheorem{proposition}[theorem]{Proposition}
\newtheorem{Examples}[theorem]{Example}
\newtheorem{corollary}[theorem]{Corollary}

\def\[{{\Big[}}
\def\]{{\Big]}}
\def\<{{\langle}}
\def\>{{\rangle}}
\def\({{\Big(}}
\def\){{\Big)}}

\def\bx{{\mathbf{x}}}

\def\dif{{\mathord{{\rm d}}}}

\def\no{\nonumber}
\def\={&\!\!=\!\!&}
\def\bt{\begin{theorem}}
\def\et{\end{theorem}}
\def\bl{\begin{lemma}}
\def\el{\end{lemma}}
\def\br{\begin{remark}}
\def\er{\end{remark}}
\def\bd{\begin{definition}}
\def\ed{\end{definition}}
\def\bp{\begin{proposition}}
\def\ep{\end{proposition}}
\def\bc{\begin{corollary}}
\def\ec{\end{corollary}}
\def\bx{\begin{Examples}}
\def\ex{\end{Examples}}

\def\mE{{\mathbb E}}
\def\mF{{\mathbb F}}

\def\mH{{\mathbb H}}

\def\mR{{\mathbb R}}

\def\sB{{\mathscr B}}

\def\sD{{\mathscr D}}

\def\sF{{\mathscr F}}

\def\geq{\geqslant}
\def\leq{\leqslant}

\allowdisplaybreaks

\begin{document}

\title{ Reflected anticipated backward stochastic differential equations  with nonlinear resistance$^*$}
\date{}
\author{ Wu Hao$^{1}$}
\thanks{$^*$ Corresponding author :  Wu Hao}

\address{Wu Hao:
School of Mathematics and Statistics,
South-Central University For Nationalities\\
Wuhan, Hubei 430000, P.R.China\\
email: wuhaomoonsky@163.com
}

\begin{abstract}
In this paper,   we consider  reflected anticipated backward stochastic differential equations (RABSDEs, for short) with an additional resistance in the generators. Firstly, we study the existence and uniqueness results. In Luo (2020), the condition of a small time horizon  is  needed. Compared with the proving method in Luo (2020), we use a different proving method to avoid requiring the Lipschitz coefficients of generators $f(t,y,z,\theta,\vartheta,m,\bar{m})$ for $y,z,\theta,\vartheta$  to be small enough.  We  only require  the Lipschitz coefficient for resistance  in generator is small enough.   Moreover,  a probabilistic structure for solution  is specified. Secondly, we give a comparison theorem for this type of equation. At last, under the linear growth condition and some other conditions on resistance , we derive the minimal solution.
\begin{flushleft}
\end{flushleft}

\begin{flushleft}
\end{flushleft}
\end{abstract}

\maketitle
\rm
\section{Introduction}
Throughout this paper,  for $x, y \in \mR,$ we use $|x |$  to denote the Euclidean norm of
$x,$ and use  $\langle x, y\rangle$ to denote the Euclidean inner product. For $B\in \mR^{ d},$     $|B |$ represents  $\sqrt{\mathrm{Tr} BB^{\ast}.}$  Let  $(\Omega, \sF, P)$ be  a complete probability space taking along
 a $d-$dimensional Brownian motion $\{W_t\}_{0\leq t\leq T}$. $\mF\doteq \{\sF_{t}\}_{t\in [0,T+K]}$ is the natural filtration generated by  $W$.
Consider the following  backward stochastic differential equations (BSDEs):
\begin{align*}
Y_t = \xi + \int^{T}_{t}f(s, Y_s, Z_s)\dif s - \int^{T}_{t}Z_s\dif W_s.
\end{align*}
In this  equation,  there exists a triple of coefficients  ($\xi$,  T,  $f$) so-called  parameters, in which $\xi$ is named  the terminal value, $T$ is a constant called the time horizon, $f$ is a random function  so-named the generator. For the case of linear type  of the above equation,  it has been introduced    by Bismut (1976), as equation for the adjoint process in the stochastic version of Pontryagin maximum principle . For the case of nonlinear type of the above equation,
  Pardaoux and Peng (1990) have established the rigorous framework for the analysis of nonlinear BSDEs. Also, Pardaoux and Peng (1990) have get
   Feynman-Kac formula which could afford a probabilistic interpretation for a classes of PDEs. Since then,  the literature on this topic bloomed, both in the direction of obtaining qualitative and quantitative results for the generalized emerging equations and on developing applications which aim to   stochastic partial differential equations,   controlled games and finance, ect. (c.f. Cr\'{e}pey and Song (2015); Madec(2015), Yang and Zhang (2014); Nunno and Sjursen (2014) ).

As the main   driving force of this article,   a new research milestone in the study of BSDEs were established by the approaches of reflected BSDEs, which were analyzed as stand alone or in connection to PDEs. For example, Karoui et al. (1997)  first introduced  Reflected BSDEs and studied  related obstacle
problems for PDEs.  Later, Ren and Otmani (2010) studied  generalized reflected BSDEs
driven by a Levy process and an obstacle problem for PDEs
with a nonlinear Neumann boundary condition.  Recently, a new class of reflected BSDEs  with nonlinear resistance has been introduced in Qian and Xu (2018), where they obtain the existence and uniqueness of  solution , which generalized the work of Karoui et al. (1997) to the case which allowed the generator to have a  resistance term.   Based on Qian and Xu (2018),   Luo (2020)  studied  reflected BSDEs with time-delayed generators and nonlinear resistance.
In addition,  much attention was paid
on studying the minimal  solutions of BSDEs (c.f.  Lepeltier and  Martin (1997), Fan et al. (2011), Fan and Jiang (2012)), whose method on studying minimal  solutions will be used in this article.

Obviously,  the generators of all the above equations didn't contain future values of solutions. In order to improving it ,  recently, Peng and Yang (2009) introduced
a new type of BSDEs, called  anticipated BSDEs (ABSDEs) as follows:
\begin{align*}
 \begin{cases}
 -\dif Y_t= f(t, Y_t, Z_t, Y_{t+\mu(t)}, Z_{t+\nu(t)})\dif t - Z_t\dif W_t, &  t\in[0, T], \\
 Y_t= \xi_t, \qquad Z_t = \eta_t, &  t\in [T, T+\delta],
 \end{cases}
\end{align*}
where $\mu (\cdot): [0, T]\rightarrow \mR^{+} \backslash \{0\}$ and $\nu(\cdot): [0, T]\rightarrow \mR^{+} \backslash \{0\}$  are continuous functions satisfying:
\begin{enumerate}[(i)]
\item \,\,   $\delta \geq 0$ is a constant such that for each $t\in [0, T],$
$$t + \mu(t) \leq T + \delta,\,\, t + \nu(t) \leq T + \delta.$$
\item \,\, There is a constant $L \geq 0$ such that for any $t\in [0, T]$ and any nonnegative integrable function $g(\cdot),$
    \begin{align*}\int^{T}_{t}g(s+ \mu(s))\dif s \leq L\int^{T+\delta}_{t}g(s)\dif s,\,\,\int^{T}_{t}g(s+ \nu(s))\dif s \leq L\int^{T+\delta}_{t}g(s)\dif s. \end{align*}
\end{enumerate}

Under the assumptions of  the Lipschitz conditions on $f$, Peng and Yang (2009) have proved  that ABSDEs  had a unique solution and given some comparison theorems. Furthermore,    by using
  the perfect  duality between SDDEs and  ABSDEs,  Peng and Yang (2009)  have solved a interesting optimal controlled problem. Since then, ABSDEs were further studied in many other articles. For examples, Lu and Ren  (2013) have studied ABSDEs driven by Markov chain. Liu and Ren (2015)  have studied anticipated BSDEs driven by time-changed L\'{e}vy noises. Feng (2016) has  studied the ABSDEs with Reflection.  Yang and   Elliot (2013) have studied the minimal solution for ABSDEs, in which the generators $f$ need to satisfied  the continuous  and linear growth  conditions.
Motivated by the above articles, we will study the following reflected anticipated backward stochastic differential equations with nonlinear resistance:
\begin{align}
 \begin{cases}
 -\dif Y_t= f(t,  Y_t, Z_t, Y_{t+\mu(t)}, Z_{t+\nu(t)}, G_{t}(K),  \mE^{\sF_{t}}[G_{t+\epsilon(t)}(K)] )\dif t+\dif K_{t} - Z_t\dif W_t,   t\in[0, T]; \\
 Y_{t}\geq S_{t}, t\in [0,T], \int^{T}_{0}(Y_{t}-S_{t})\dif K_{t}=0;\\
 Y_t= \xi_t,  Z_t = \eta_t, t\in [T, T+\delta];K_t = \zeta_t , t\in (T, T+\delta].
 \end{cases}
\end{align}
  The resistance term is allowed to depend on the past values and future values of the increasing process.

We close this part by giving our aims in this article.  The first  aim in this paper is to get the  existence and uniqueness results of adapted solution as well as  a probabilistic structure  is specified for solution.

The second aim is to derive a comparison theorem.

The last aim is to get  the minimal solution of the above equation with linear growth and continuous condition.

\section{Preliminaries}
\subsection{Notations }
For Euclidean space $\mH,$ we introduce the following spaces:

$L^{2}(\sF_{t}; \mH)$ is represented as a  space of $\mH$-valued $\sF_{t}-$measurable random variables $\phi$ satisfying $\|\phi\|_{2}\doteq (\mE[|\phi|^{2}])^{\frac{1}{2}}< \infty.$\\

$L^{2}_{\mF}(0, T; \mH)$ is represented as a  space of $\mH$-valued $\mF-$adapted stochastic processes $\{\varphi_{s}, s\in [0, T]\}$ satisfying $\|\varphi\|_{L^{2}_{\mF}(0, T)}\doteq \bigg(\bigg. \mE\bigg[\bigg. \int^{T}_{0}|\varphi_s|^{2}\dif s\bigg]\bigg. \bigg)\bigg.^{\frac{1}{2}}< \infty.$\\

$S^{2}_{\mF}(0, T; \mH)$ is represented as a  space of continuous  processes $\{\varphi_{s}, s\in [0, T]\}$ in $L^{2}_{\mF}(0, T; \mH)$ satisfying $\|\varphi\|_{S^{2}_{\mF}(0, T)}\doteq \bigg(\bigg. \mE [\sup_{0 \leqslant s \leqslant T}|\varphi_s|^{2}]\bigg)\bigg.^{\frac{1}{2}}< \infty.$\\

$H^{2}_{\mF}(0, T; \mH)$ is represented as a  space of $\mH$-valued $\mF-$adapted stochastic processes $\{\varphi_{s}, s\in [0, T]\}$ in $L^{2}_{\mF}(0, T; \mH)$ satisfying $\|\varphi\|_{H^{2}_{\mF}(0, T)}\doteq \bigg(\bigg. \mE [\sup_{0 \leqslant s \leqslant T}|\varphi_s|^{2}]\bigg)\bigg.^{\frac{1}{2}}< \infty.$\\

$\sD(\mR)$ be denoted by all the  functions with only a finite number of discontinuities from $[0,T+K]$ to $\mR.$

$G_{t}$ is a function from $\sD(\mR)$ to $\mR$ for any $t\in [0, T+K].$

$
\mbox{Let} \,\, \beta>0\,\, \mbox{be a constant},\,\, L^{2}_{\mF}(0, T; \mH; \beta)$ is represented as a  space of $\mH$-valued $\mF-$adapted stochastic processes $\{\varphi_{s}, s\in [0, T]\}$ satisfying $\|\varphi\|_{L^{2}_{\mF}(0, T;  \beta)}\doteq \bigg(\bigg. \mE\bigg[\bigg. \int^{T}_{0}e^{\beta s}|\varphi_s|^{2}\dif s\bigg]\bigg. \bigg)\bigg.^{\frac{1}{2}}< \infty.$\\

$
\mbox{Let} \,\, \beta>0\,\, \mbox{be a constant},\,\, S^{2}_{\mF}(0, T; \mH; \beta)$  is represented as a  space of continuous  processes $\{\varphi_{s}, s\in [0, T]\}$ in $L^{2}_{\mF}(0, T; \mH; \beta)$  satisfying $\|\varphi\|_{S^{2}_{\mF}(0, T;  \beta)}\doteq \bigg(\bigg. \mE [\sup_{0 \leqslant s \leqslant T}e^{\beta s}|\varphi_s|^{2}]\bigg)\bigg.^{\frac{1}{2}}< \infty.$\\

$
\mbox{Let} \,\, \gamma>0, \beta>0\,\, \mbox{be two constants},\,\, H^{2}_{\mF}(0, T; \mH;\gamma, \beta)$ is represented as a  space of $\mH$-valued $\mF-$adapted stochastic processes $\{\varphi_{s}, s\in [0, T]\}$ in $L^{2}_{\mF}(0, T; \mH)$ satisfying $\|\varphi\|_{H^{2}_{\mF}(0, T;\gamma, \beta)}\doteq \bigg(\bigg. \mE [\frac{1}{\gamma}\sup_{0 \leqslant s \leqslant T}e^{\beta s}|\varphi_s|^{2}]\bigg)\bigg.^{\frac{1}{2}}< \infty.$\\

Obviously, $\|\cdot\|_{H^{2}_{\mF}(0, T;\gamma, \beta)}$  and $\|\cdot\|_{H^{2}_{\mF}(0, T)}$ are two equivalent norms. $\|\cdot\|_{S^{2}_{\mF}(0, T; \beta)}$  and $\|\cdot\|_{S^{2}_{\mF}(0, T)}$ are two equivalent norms. $\|\cdot\|_{L^{2}_{\mF}(0, T)}$ and $\|\cdot\|_{L^{2}_{\mF}(0, T;  \beta)}$ are two equivalent norms. In addition, let $\mE^{\sF_{t}}[\cdot]\doteq \mE[\cdot|\sF_{t}].$

\subsection{Hypotheses}
\begin{enumerate}[(H1)]
\item Assume that for any
 $t\in [0, T],$ $f(t,\omega, y ,z, \theta, \vartheta,m,\bar{m}): \Omega\times \mR\times \mR^{d} \times L^{2}( \sF_{r} ; \mR)\times L^{2}(\sF_{\bar{r}}; \mR^{d})\times \mR \times \mR \rightarrow L^{2}(\sF_{t}; \mR),$ where
$r, \bar{r} \in [t, T+\delta],$ and $f$ satisfies the following conditions:
$$\mE\left[\int^{T}_{0}|f(s,0,0,0 ,0, 0,0)|^{2}\dif s \right]< \infty,$$
\item For any $t\in [0, T],  y,y^{\prime},m,m^{\prime},\bar{m},\bar{m}^{\prime}\in \mR, z, z^{\prime}\in \mR^{d}, \theta, \theta^{\prime}\in L^{2}( \sF_{r} ; \mR),\vartheta, \vartheta^{\prime}\in L^{2}(\sF_{\bar{r}}; \mR^{d}), r, \bar{r} \in [t, T + K],$ we have
    \begin{align*}
    |&f(t, y,z,\theta, \vartheta,m, \bar{m }) - f(t, y^{\prime},z^{\prime},m^{\prime},\theta^{\prime}, \vartheta^{\prime},m', \bar{m}')|\\
    &\leq C(|y- y^{\prime}| + |z-z^{\prime}| + \mE^{\sF_{t}}[|\theta - \theta^{\prime}|+|\vartheta - \vartheta^{\prime}|])+C_{1}(|m- m^{\prime}|+|\bar{m}- \bar{m}^{\prime}|).
    \end{align*}
\item Assume that $\xi, \zeta\in S^{2}_{\mF}(T, T + \delta; \mR) $ and $\eta\in L^{2}_{\mF}(T, T + \delta; \mR), S \in S^{2}_{\mF}(0, T; \mR) .$

\item For any $ y,y^{\prime}\in \sD(\mR),$ $G_{t}$ satisfies $G(0)=0$  and
$$G_{t}(y)=G_{t}(\{y_{s\wedge t}\}_{0\leq s \leq T}), $$
$$|G_{t}(y)-G_{t}(y')|\leq \sup_{0\leq s\leq t}|y_{s}-y'_{s}|.$$
    \end{enumerate}

\br
Given some $\varepsilon \geq 0,$ we give some examples satisfying $(\mathrm{H}4):$

$$ 1^{\circ} G_{t}: y\rightarrow \int^{t}_{\frac{t}{2}}y^{+}_{s}\dif s,2^{\circ} G_{t}: y\rightarrow \sup_{ \frac{t}{2}\leq u\leq t}y_{s},3^{\circ}G_{t}: y\rightarrow y^{+}_{(t-\varepsilon)^{+}},4^{\circ} G_{t}: y\rightarrow \int^{t}_{\frac{t}{2}}y^{+}_{s}\dif s, 5^{\circ} G_{t}: y\rightarrow \frac{1}{\varepsilon}\int^{t}_{(t-\varepsilon)^{+}}y^{+}_{s}\dif s.$$

\er
\bd
We say that $(Y, Z, K)$ is a solution of Eq.(1), if the following conditions hold:
\begin{enumerate}[(a)]
\item (Y, Z, K) satisfies  the Eq.(1)
 \item $(Y, Z, K)\in S^{2}_{\mF}(0, T+\delta; \mR)\times L^{2}_{\mF}(0, T+\delta; \mR^{d}) \times H^{2}_{\mF}(0, T+\delta; \mR). $
 \item $K_{t} , t\in [0,T]$  is a continuous increasing process.

    \end{enumerate}

\ed

We close this preliminaries by introducing a convention.  A  convention is needed as follows:  $C^{\prime}$  is a positive constant and its value could be allowed to vary from one place to another but $C^{\prime}$ only rely on the constants in the assumptions.

\section{main results}
\subsection{Existence and uniqueness result}In this part, we will study the existence and uniqueness result for Eq.(1).

\bt Let $\lambda =4(6C^{2}+6C^{2}L), \beta=\lambda +2, \gamma= 4 [(4Te^{\beta T}+16Te^{\beta T}(1+e^{\beta T})) (6C^{2}+6C^{2}L)+4e^{\beta T}].$  $C_{1}$ is
small enough such that $6\bigg(\bigg.4Te^{\beta T}+16Te^{\beta T}(1+e^{\beta T})+\frac{\gamma}{\lambda}\bigg)\bigg.(C^{2}_{1}T+C^{2}_{1}(T+\delta)L)\leq \frac{1}{4}.$
  Under assumptions $(\mathrm{H}1)-(\mathrm{H}4), (\mathrm{i}),(\mathrm{ii}) $, there exists a unique triple of sulution  $(Y, Z, K)\in S^{2}_{\mF}(0, T+\delta; \mR)\times L^{2}_{\mF}(0, T+\delta; \mR^{d}) \times H^{2}_{\mF}(0, T+\delta; \mR)$ satisfying Eq.(1) and $Y$ has the following representation.
 \begin{align*}
Y_{t}=\mathrm{ess}\sup_{\tau \in \Gamma_{t}}\mE\bigg[\bigg. \int^{\tau}_{t} f(s,  Y_s, Z_s, Y_{s+\mu(s)}, Z_{s+\nu(s)}, G_{s}(K), \mE^{\sF_{s}}[G_{s+\epsilon(s)}(K)])\dif s + S_{\tau}1_{\tau<T} +\xi_{T}1_{\tau=T}\bigg| \bigg. \sF_{t}   \bigg]\bigg., t\in [0,T],
\end{align*}
where $\Gamma$ is the set of all stopping times taking values in $[0,T]$  and $\Gamma_{t}=\{\tau \in \Gamma: \tau\geq t\}.$
\et
\begin{proof}
Set $\sB^{2}_{\upsilon}=\{(U,V, k)|(U,V,k)\in L^{2}_{\mF}(0, T+\delta; \mR;\beta)\times L^{2}_{\mF}(0, T+\delta; \mR^{d};\beta)\times H^{2}_{\mF}(0, T+\delta; \mR;\gamma, \beta); U_s=\xi_{s}, V_{s}=\eta_{s}, s\in [T, T+K], k_t = \zeta_t, s\in (T, T+K]\}.$ For any $(U,V,k)\in \sB^{2}_{\upsilon},$  by Theorem 5.2 in Karoui et al. (1997),  the following equation has a unique solution $(Y,Z,K).$
\begin{align}
 \begin{cases}
  Y_t=\int^{T}_{t} f(s,  Y_s, Z_s, Y_{s+\mu(s)}, Z_{s+\nu(s)}, G_{s}(k), \mE^{\sF_{s}}[G_{s+\epsilon(s)}(k)])\dif s+ K_{T}-K_{t} - \int^{T}_{t}Z_s\dif W_s,   t\in[0, T]; \\
 Y_{s}\geq S_{s}, t\in [0,T], \int^{T}_{0}(Y_{s}-S_{s})\dif K_{s}=0;\\
 Y_t= \xi_t,  Z_t = \eta_t, t\in [T, T+\delta]; K_t = \zeta_t, t\in (T, T+\delta],
 \end{cases}
\end{align}
and
\begin{align*}
Y_{t}=\mathrm{ess}\sup_{\tau \in \Gamma_{t}}\mE\bigg[\bigg. \int^{\tau}_{t} f(s,  Y_s, Z_s, Y_{s+\mu(s)}, Z_{s+\nu(s)}, G_{s}(k), \mE^{\sF_{s}}[G_{s+\epsilon(s)}(k)])\dif s + S_{\tau}1_{\tau<T} +\xi_{T}1_{\tau=T}\bigg| \bigg. \sF_{t}   \bigg]\bigg., t\in [0,T].
\end{align*}
We will use the Fixed Point Theorem to prove the above theorem.  Define a map $F$ from $\sB^{2}_{\upsilon}$ to $\sB^{2}_{\upsilon}$ as follows:
 for $(U,V,k)\in \sB^{2}_{\upsilon},$  $(Y,Z,K)= F(U,V,k).$ For any $(U,V,k),(U',V',k')\in \sB^{2}_{\upsilon}, $  define $(Y,Z,K)= F(U,V,k), (Y',Z',K')= F(U',V',k').$  Furthermore,  $Y$ and $Y'$ have the following representations:
   \begin{align*}
Y_{t}=\mathrm{ess}\sup_{\tau \in \Gamma_{t}}\mE\bigg[\bigg. \int^{\tau}_{t} f(s,  Y_s, Z_s, Y_{s+\mu(s)}, Z_{s+\nu(s)}, G_{s}(k), \mE^{\sF_{s}}[G_{s+\epsilon(s)}(k)])\dif s + S_{\tau}1_{\tau<T} +\xi_{T}1_{\tau=T}\bigg| \bigg. \sF_{t}   \bigg]\bigg., t\in [0,T],
\end{align*}
  and
  \begin{align*}
Y'_{t}=\mathrm{ess}\sup_{\tau \in \Gamma_{t}}\mE\bigg[\bigg. \int^{\tau}_{t} f(s,  Y'_s, Z'_s, Y'_{s+\mu(s)}, Z'_{s+\nu(s)}, G_{s}(k'), \mE^{\sF_{s}}[G_{s+\epsilon(s)}(k')])\dif s + S_{\tau}1_{\tau<T} +\xi_{T}1_{\tau=T}\bigg| \bigg. \sF_{t}   \bigg]\bigg., t\in [0,T].
\end{align*}
  Set
 $$\hat{U}= U-U',\hat{V}= V-V', \hat{Y}= Y-Y', \hat{Z}= Z-Z',\hat{K}= K-K', \hat{k}= k-k',$$
 $$\hat{f}_s= f(s,  Y_s, Z_s, Y_{s+\mu(s)}, Z_{s+\nu(s)}, G_{s}(k), \mE^{\sF_{s}}[G_{s+\epsilon(s)}(k)])- f(s,  Y'_s, Z'_s, Y'_{s+\mu(s)}, Z'_{s+\nu(s)}, G_{s}(k'), \mE^{\sF_{s}}[G_{s+\epsilon(s)}(k')]). $$
 Consider the following equations:
\begin{align*}
 \hat{Y}_t=\int^{T}_{t} \hat{f}_s\dif s + \hat{K}_{T}-\hat{K}_{t}- \int^{T}_{t}\hat{Z}_s\dif W_s,   t\in[0, T].
 \end{align*}
Using It\^{o}'s formula to $e^{\beta t}|Y_{t}|^{2},$ we have
\begin{align*}
 \mE&[e^{\beta t}|\hat{Y}_{t}|^{2}] + \beta\mE\bigg[\bigg.\int^{T}_{t}e^{\beta s}|\hat{Y}_{s}|^{2}\dif s\bigg]\bigg. +\mE\bigg[\bigg.\int^{T}_{t}e^{\beta s}|\hat{Z}_{s}|^{2}\dif s \bigg]\bigg.
 =  2\mE\bigg[\bigg.\int^{T}_{t}e^{\beta s}\langle \hat{Y}_{s},  \hat{f}_s\rangle\dif s\bigg]\bigg..
\end{align*}
Then,
\begin{align}
 \mE[e^{\beta t}|\hat{Y}_{t}|^{2}] + \beta\mE\bigg[\bigg.\int^{T}_{t}e^{\beta s}|\hat{Y}_{s}|^{2}\dif s\bigg]\bigg. +\mE\bigg[\bigg.\int^{T}_{t}e^{\beta s}|\hat{Z}_{s}|^{2}\dif s \bigg]\bigg.
  &\leq \lambda\mE\bigg[\bigg.\int^{T}_{t}e^{\beta s}| \hat{Y}_{s}|^{2}\dif s\bigg]\bigg.+\frac{1}{\lambda}\mE\bigg[\bigg.\int^{T}_{t}e^{\beta s}  |\hat{f}_s|^{2}\dif s\bigg]\bigg.
\end{align}
Since,
\begin{align*}
\hat{K}_{t}=-\hat{Y}_{t}+\hat{Y}_{0}-\int^{t}_{0}\hat{f}_{s}\dif s + \int^{t}_{0}\hat{Z}_{s}\dif W_{s}.
\end{align*}
Thus,
\begin{align*}
\frac{1}{\gamma}e^{\beta t}|\hat{K}_{t}|^{2}\leq\frac{4}{\gamma}e^{\beta t}|\hat{Y}_{t}|^{2}+\frac{4}{\gamma}e^{\beta t}|\hat{Y}_{0}|^{2}+\frac{4T}{\gamma}e^{\beta t}\int^{t}_{0}|\hat{f}_{s}|^{2}\dif s + \frac{4}{\gamma}e^{\beta t}\bigg|\bigg. \int^{t}_{0}\hat{Z}_{s}\dif W_{s}\bigg|\bigg.^{2}.
\end{align*}
Using BDG's inequality, we derive
\begin{align}
\mE\bigg[\bigg.\sup_{0\leq t\leq T}\frac{1}{\gamma}e^{\beta t}|\hat{K}_{t}|^{2}\bigg]\bigg.\leq  4(1+e^{\beta T})\mE\bigg[\bigg.\sup_{0\leq t\leq T}\frac{1}{\gamma}e^{\beta t}|\hat{Y}_{t}|^{2}\bigg]\bigg.+\frac{4Te^{\beta T}}{\gamma}\mE\bigg[\bigg.\int^{T}_{0}e^{\beta s}|\hat{f}_{s}|^{2}\dif s\bigg]\bigg. + \frac{4e^{\beta T}}{\gamma}\mE\bigg[\bigg. \int^{T}_{0}e^{\beta s}|\hat{Z}_{s}|^{2}\dif s\bigg]\bigg..
\end{align}
From the following representations:
 \begin{align*}
Y_{t}=\mathrm{ess}\sup_{\tau \in \Gamma_{t}}\mE\bigg[\bigg. \int^{\tau}_{t} f(s,  Y_s, Z_s, Y_{s+\mu(s)}, Z_{s+\nu(s)}, G_{s}(k), \mE^{\sF_{s}}[G_{s+\epsilon(s)}(k)])\dif s + S_{\tau}1_{\tau<T} +\xi_{T}1_{\tau=T}\bigg| \bigg. \sF_{t}   \bigg]\bigg., t\in [0,T],
\end{align*}
and
\begin{align*}
Y_{t}=\mathrm{ess}\sup_{\tau \in \Gamma_{t}}\mE\bigg[\bigg. \int^{\tau}_{t} f(s,  Y_s, Z_s, Y_{s+\mu(s)}, Z_{s+\nu(s)}, G_{s}(k), \mE^{\sF_{s}}[G_{s+\epsilon(s)}(k)])\dif s + S_{\tau}1_{\tau<T} +\xi_{T}1_{\tau=T}\bigg| \bigg. \sF_{t}   \bigg]\bigg., t\in [0,T],
\end{align*}
we deduce

\begin{align*}
|\hat{Y}_{t}|\leq \mathrm{ess}\sup_{\tau \in \Gamma_{t}}\mE\bigg[\bigg. \int^{\tau}_{t} |\hat{f}_s|\dif s \bigg| \bigg. \sF_{t}   \bigg]\bigg..
\end{align*}
Thus,
\begin{align*}
|\hat{Y}_{t}|
\leq \bigg[\bigg. \int^{T}_{0} |\hat{f}_s|\dif s \bigg| \bigg. \sF_{t}   \bigg]\bigg..
\end{align*}
It leads to
\begin{align*}
|\hat{Y}_{t}|^{2}
 \leq \bigg(\bigg. \int^{T}_{0} |\hat{f}_s|\dif s \bigg| \bigg. \sF_{t}   \bigg)\bigg.^{2}.
\end{align*}
Doob's maximal inequality implies that
\begin{align*}
\mE\bigg[\bigg.\sup_{0\leq t \leq T}|\hat{Y}_t|^{2}\bigg]\bigg. \leq  4T\mE\bigg[\bigg.\int^{T}_{0} |\hat{f}_s|^{2}\dif s\bigg]\bigg..
 \end{align*}
Then, we have

 \begin{align*}
\mE\bigg[\bigg.\sup_{0\leq t \leq T}\frac{1}{\gamma}e^{\beta t}|\hat{Y}_t|^{2}\bigg]\bigg. \leq  \frac{4Te^{\beta T}}{\gamma}\mE\bigg[\bigg.\int^{T}_{0} e^{\beta s}|\hat{f}_s|^{2}\dif s\bigg]\bigg..
 \end{align*}
By (4),  it holds that
\begin{align}
\mE\bigg[\bigg.\sup_{0\leq s\leq T}\frac{1}{\gamma}e^{\beta s}|\hat{K}_{s}|^{2}\bigg]\bigg.&\leq  (1+e^{\beta T})\frac{16Te^{\beta T}}{\gamma}\mE\bigg[\bigg.\int^{T}_{0} e^{\beta s}|\hat{f}_s|^{2}\dif s\mE\bigg]\bigg.  +\frac{4Te^{\beta T}}{\gamma}\mE\bigg[\bigg.\int^{T}_{0}e^{\beta s}|\hat{f}_{s}|^{2}\dif s\bigg]\bigg.
  + \frac{4e^{\beta T}}{\gamma}\mE\bigg[\bigg. \int^{T}_{0}|e^{\beta s}\hat{Z}_{s}|\dif s\bigg]\bigg..\no\\
&\leq \frac{4Te^{\beta T}+16Te^{\beta T}(1+e^{\beta T})}{\gamma}\mE\bigg[\bigg.\int^{T}_{0}e^{\beta s}|\hat{f}_{s}|^{2}\dif s\bigg]\bigg. + \frac{4e^{\beta T}}{\gamma}\mE\bigg[\bigg. \int^{T}_{0}e^{\beta s}|\hat{Z}_{s}|^{2}\dif s\bigg]\bigg..
\end{align}
Combining  (3) and (5), one gets that
 \begin{align*}
&\mE\bigg[\bigg.\sup_{0\leq s\leq T}\frac{1}{\gamma}e^{\beta s}|\hat{K}_{s}|^{2}\bigg]\bigg.+\beta\mE\bigg[\bigg.\int^{T}_{0}e^{\beta s}|\hat{Y}_{s}|^{2}\dif s\bigg]\bigg.+\mE\bigg[\bigg.\int^{T}_{0}e^{\beta s}|\hat{Z}_{s}|^{2}\dif s\bigg]\bigg.\\
&\leq\lambda\mE\bigg[\bigg.\int^{T}_{0}e^{\beta s}| \hat{Y}_{s}|^{2}\dif s\bigg]\bigg.
+\bigg(\bigg.\frac{4Te^{\beta T}+16Te^{\beta T}(1+e^{\beta T})}{\gamma}+\frac{1}{\lambda}\bigg)\bigg.\mE\bigg[\bigg.\int^{T}_{0}e^{\beta s}|\hat{f}_{s}|^{2}\dif s\bigg]\bigg. + \frac{4e^{\beta T}}{\gamma}\mE\bigg[\bigg. \int^{T}_{0}e^{\beta s}|\hat{Z}_{s}|^{2}\dif s\bigg]\bigg..
\end{align*}
Hence,
\begin{align*}
&\mE\bigg[\bigg.\sup_{0\leq s\leq T}\frac{1}{\gamma}e^{\beta s}|\hat{K}_{s}|^{2}\bigg]\bigg.+\beta\mE\bigg[\bigg.\int^{T}_{0}e^{\beta s}|\hat{Y}_{s}|^{2}\dif s\bigg]\bigg.+\mE\bigg[\bigg.\int^{T}_{0}e^{\beta s}|\hat{Z}_{s}|^{2}\dif s\bigg]\bigg.\\
&\leq\bigg(\bigg. \lambda + \bigg(\bigg.\frac{4Te^{\beta T}+16Te^{\beta T}(1+e^{\beta T})}{\gamma}+\frac{1}{\lambda}\bigg)\bigg.(6C^{2}+6C^{2}L) \bigg)\bigg.\mE\bigg[\bigg.\int^{T}_{0}e^{\beta s}| \hat{Y}_{s}|^{2}\dif s\bigg]\bigg.\\
& \quad  + \bigg(\bigg.\frac{4e^{\beta T}}{\gamma} + \bigg(\bigg.\frac{4Te^{\beta T}+16Te^{\beta T}(1+e^{\beta T})}{\gamma}+\frac{1}{\lambda}\bigg)\bigg.(6C^{2}+6C^{2}L) \bigg)\bigg.   \mE\bigg[\bigg. \int^{T}_{0}e^{\beta s}|\hat{Z}_{s}|^{2}\dif s\bigg]\bigg.\\
& \quad + 6\bigg(\bigg.4Te^{\beta T}+16Te^{\beta T}(1+e^{\beta T})+\frac{\gamma}{\lambda}\bigg)\bigg.(C^{2}_{1}T+C^{2}_{1}(T+\delta)L)\mE\bigg[\bigg.\sup_{0\leq s\leq T}\frac{1}{\gamma}e^{\beta s}|\hat{k}_{s}|^{2}\bigg]\bigg..
\end{align*}
Choosing $\lambda =4(6C^{2}+6C^{2}L), \beta=\lambda +2, \gamma= 4 [(4Te^{\beta T}+16Te^{\beta T}(1+e^{\beta T})) (6C^{2}+6C^{2}L)+4e^{\beta T}],$  we have
$$\mE\bigg[\bigg.\sup_{0\leq s\leq T}\frac{1}{\gamma}e^{\beta s}|\hat{K}_{s}|^{2}\bigg]\bigg.+\mE \bigg[\bigg.\int^{T}_{0}e^{\gamma s}|\hat{Y}_{s}|^{2}\dif s\bigg]\bigg. +\mE\bigg[\bigg.\int^{T}_{0}e^{\gamma s}|\hat{Z}_{s}|^{2}\dif s \bigg]\bigg.\leq \frac{1}{2}\mE\bigg[\bigg.\sup_{0\leq s\leq T}\frac{1}{\gamma}e^{\beta s}|\hat{k}_{s}|^{2}\bigg]\bigg.. $$
We know that $F$ is a strict contraction mapping on $\sB^{2}_{\upsilon}.$  Thus, there exists a unique pair $(Y, Z,K)\in L^{2}_{\mF}(0, T+\delta; \mR)\times L^{2}_{\mF}(0, T+\delta; \mR^{d})\times H^{2}_{\mF}(0, T+\delta; \mR) $ satisfying Eq.(1) and Using BDG's inequality, we get $(Y, Z, K)\in S^{2}_{\mF}(0, T+\delta; \mR)\times L^{2}_{\mF}(0, T+\delta; \mR^{d})\times H^{2}_{\mF}(0, T+\delta; \mR).$ The proof   is complete.
\end{proof}
\br
In Luo (2020), under their method, the condition of a small time horizon is  needed. In essence, the Lipschitz coefficients of generator $f(t,y,z,\theta, \vartheta,m,\bar{m})$ for $y, z, \theta, \vartheta,m,\bar{m}$ are required to be small enough. Compared with the proving method in Luo (2020), we introduce three useful constants $\beta,\gamma,\lambda$ and a space  $H^{2}_{\mF}(0, T; \mH;\gamma, \beta)$  in the proof and use a different proving method to avoid requiring the Lipschitz coefficients for $y, z, \theta,\vartheta$ to be small enough.     We  only require  the Lipschitz coefficient for resistance  in generator is small enough.  This is the difference.
\er

\subsection{Comparison theorem}Next, we will give a comparison theorem for the following equations.
Let $(Y^{(1)}, Z^{(1)}, K^{(1)})$ and $(Y^{(2)}, Z^{(2)}, K^{(2)})$ be the solutions of the following two 1-dimensional RABSDEs, respectively:
 \begin{align*}
 \begin{cases}
 -\dif Y^{(j)}_{t}= f_{j}(t, Y^{(j)}_{t}, Z^{(j)}_{t}, Y^{(j)}_{t+\mu(t)}, G_{t}(K^{(j)}), \mE^{\sF_{t}}[G_{t+\epsilon(t)}(K^{(j)})])\dif t +\dif K^{(j)}_{t}- Z^{(j)}_{t}\dif W_{t},  t\in[0, T], \\
  Y^{(j)}_{t}\geq S_{t}, t\in [0,T], \int^{T}_{0}(Y^{(j)}_{t}-S_{t})\dif K^{(j)}_{t}=0;\\
 Y^{(j)}_{t}= \xi^{(j)}_{t}, t\in [T, T+\delta]; K^{(j)}_{t}= \zeta^{(j)}_{t}, t\in (T, T+\delta],
\end{cases}
\end{align*}
where $j= 1,2.$ For simplicity, we rewrite the above equations as the follows:
 \begin{align*}
 \begin{cases}
 -\dif Y^{(j)}_{t}= f_{j}(t, Y^{(j)}_{t}, Z^{(j)}_{t}, Y^{(j)}_{t+\mu(t)}, G_{t}(K^{(j)}), \mE^{\sF_{t}}[G_{t+\epsilon(t)}(K^{(j)})])\dif t +\dif K^{(j)}_{t}- Z^{(j)}_{t}\dif W_{t},t\in[0, T], \\
  Y^{(j)}_{t}= \xi^{(j)}_{t}, t\in [T, T+\delta]; K^{(j)}_{t}= \zeta^{(j)}_{t}, t\in (T, T+\delta].\\
\end{cases}
\end{align*}
Below, we give a comparison theorem for RABSDEs and make the following assumptions.
\begin{enumerate}[(b1)]
\item For any $t\in [0, T],y, m, \bar{m}\in \mR, z\in \mR^{d},$ $f(t,y,z,\cdot,m, \bar{m})$ is increasing, that is, $f(t,y,z,\theta_r,m, \bar{m})\geq f(t,y,z,\theta^{\prime}_{r},m, \bar{m}),$ if $\theta_r \geq \theta^{\prime}_{r},$ $\theta, \theta^{\prime}\in L_{\mF}^{2}(t, T+K; \mR),\,\,r \in [t, T+K].$
\item  For any $t\in [0, T], \theta \in L_{\mF}^{2}(t, T+K; \mR), r \in [t, T+K] ,\bar{m},y\in \mR ,z\in \mR^{d},$ $f(t,y,z,\theta_r,\cdot, \bar{m})$ is decreasing, that is $f(t,y,z,\theta_r, m,\bar{m}) \leq f(t,y,z, \theta_r, m', \bar{m}),$ if  $m \geq m^{\prime}, m , m^{\prime}\in \mR.$

\item  For any $t\in [0, T], \theta \in L_{\mF}^{2}(t, T+K; \mR), r \in [t, T+K] ,m ,y\in \mR ,z\in \mR^{d},$ $f(t,y,z,\theta_r, m,\cdot)$ is decreasing, that is $f(t,y,z,\theta_r, m,\bar{m}) \leq f(t,y,z, \theta_r, m,\bar{ m}'),$ if  $\bar{m} \geq \bar{m}^{\prime}, \bar{m} , \bar{m}^{\prime}\in \mR.$

\item If $y^{(1)}, y^{(2)}\in \sD(\mR)$ and $y^{(1)}_{t}\geq y^{(2)}_{t},$  we have $G_{t}(y^{(1)})\geq G_{t}(y^{(2)}).$
\end{enumerate}
\bt
Assume that $(f_{1}, G,\xi^{(1)}, \zeta^{(1)})$ (resp. $(f_{2}, G, \xi^{(2)}, \zeta^{(2)})$) satisfies $(\mathrm{H}1)-(\mathrm{H}4)$ (resp. $(\mathrm{H}1)-(\mathrm{H}4),$ $(\mathrm{b}1)-(\mathrm{b}4)$ )  and
$\mu$ satisfies $(\mathrm{i}), (\mathrm{ii}).$   If $\xi^{(1)}_t\geq \xi^{(2)}_t,\zeta^{(1)}_t\leq\zeta^{(2)}_t,t\in [T, T+\delta]$ and $f_{1}(t,y,z,\theta_r,m,\bar{m})\geq f_{2}(t,y,z,\theta_r,m,\bar{m})$ for any $t\in [0,T], y,m,\bar{m}\in \mR, z\in\mR^{d}, \theta \in L^{2}_{\mF}(t, T+\delta; \mR), r\in [t, T+\delta],$ then
$$Y^{(1)}_{t} \geq Y_{t}^{(2)},K^{(1)}_{t} \leq K_{t}^{(2)} ,\,a.e.,\,\, a.s.$$
\et
\begin{proof}
From $(\mathrm{b}1)-(\mathrm{b}4),$
for any  $t\in [0, T], y\in\mR, z\in \mR^{ d},$ we get  that $f_{2}(t, y, z, \theta_r ,m,\bar{m} )\geq f_{2}(t, y, z,\theta^{\prime}_r, m',\bar{m}'),$ if $ m \leq m', \theta_r \geq \theta^{\prime}_{r},\bar{m} \leq \bar{m}',$ $m,m',\bar{m},\bar{m}'\in \mR, \theta, \theta^{\prime}\in L_{\mF}^{2}(t, T+\delta; \mR),\,\,r \in [t, T+\delta].$ \\
Set
\begin{align*}
 \begin{cases}
  Y^{(3)}_{t}= \xi^{(2)}_T+\int^{T}_{t}f_{2}(s, Y^{(3)}_{s}, Z^{(3)}_{s}, Y^{(1)}_{s+\mu(s)}), G_{s}(K^{(1)}), \mE^{\sF_{s}}[G_{s+\epsilon(s)}(K^{(1)})])\dif s + K^{(3)}_{T}- K^{(3)}_{t}- \int^{T}_{t}Z^{(3)}_s\dif W_s,   t\in[0, T], \\
Y^{(3)}_{t}= \xi^{(2)}_{t}, t\in [T, T+\delta]; K^{(3)}_{t}= \zeta^{(2)}_{t}, t\in (T, T+\delta].\\
  \end{cases}
\end{align*}
From the proof of Theorem 3.1,  there exists a unique pair of $\mF-$adapted processes $(Y^{(3)}, Z^{(3)},K^{(3)} )$ $\in S_{\mF}^{2}(0, T+\delta; \mR)\times  L_{\mF}^{2}(0, T; \mR^{d})\times H_{\mF}^{2}(0, T+\delta; \mR)$ that satisfies the above BSDE. Since $$f_{1}(s,y,z,Y^{(1)}_{s + \mu (s)}, G_{s}(K^{(1)}), \mE^{\sF_{s}}[G_{s+\epsilon(s)}(K^{(1)})])\geq f_{2}(s,y,z,Y^{(1)}_{s + \mu (s)}, G_{s}(K^{(1)}), \mE^{\sF_{s}}[G_{s+\epsilon(s)}(K^{(1)})]), $$ from Theorem 4.1 and Theorem 5.2 in Karoui et al. (1997),  we have

$$Y^{(1)}_t\geq Y^{(3)}_t, K^{(1)}_t\leq K^{(3)}_t,\,\,a.e., \,\,a.s.$$
Set\\
\begin{align*}
 \begin{cases}
  Y^{(4)}_t= \xi^{(2)}_T+\int^{T}_{t}f_{2}(s, Y^{(4)}_{s}, Z^{(4)}_{s}, Y^{(3)}_{s+\mu(s)},G_{s}(K^{(3)}), \mE^{\sF_{s}}[G_{s+\epsilon(s)}(K^{(3)})])\dif s + K^{(4)}_{T}- K^{(4)}_{t}- \int^{T}_{t}Z^{(4)}_s\dif W_s, &  t\in[0, T], \\
 Y^{(4)}_{t}= \xi^{(2)}_{t}, t\in [T, T+\delta]; K^{(4)}_{t}= \zeta^{(2)}_{t}, t\in (T, T+\delta].\\
  \end{cases}
\end{align*}
Since $Y^{(1)}_t\geq Y^{(3)}_t, K^{(1)}_t\leq K^{(3)}_t,$  by  $(\mathrm{b}1)-(\mathrm{b}4)$ and comparison  theorem  in Karoui et al. (1997), we know that
$$Y^{(3)}_t\geq Y^{(4)}_t, K^{(3)}_t\leq K^{(4)}_t,\,\,a.e., \,\,a.s.$$
For $n=5,6,7,\ldots,$ we consider the following classical BSDE:
\begin{align*}
 \begin{cases}
  Y^{(n)}_t= \xi^{(2)}_T+\int^{T}_{t}f_{2}(s, Y^{(n)}_{s}, Z^{(n)}_{s}, Y^{(n-1)}_{s+\mu(s)},G_{s}(K^{(n-1)}), \mE^{\sF_{s}}[G_{s+\epsilon(s)}(K^{(n-1)})])\dif s + K^{(n)}_{T}- K^{(n)}_{t}- \int^{T}_{t}Z^{(n)}_s\dif W_s,  t\in[0, T], \\
Y^{(n)}_{t}= \xi^{(2)}_{t}, t\in [T, T+\delta]; K^{(n)}_{t}= \zeta^{(2)}_{t}, t\in (T, T+\delta].\\
  \end{cases}
\end{align*}
Similarly,
$$Y^{(4)}_t \geq Y^{(5)}_t\geq\cdots \geq Y^{(n)}_t \geq ,\cdots,\,\,a.e.,\,\,a.s.$$
$$K^{(4)}_t \leq K^{(5)}_t\geq\cdots \leq K^{(n)}_t \leq,\cdots,\,\,a.e.,\,\,a.s.$$
Using the similar method in the proof of Theorem 3.1, we know that $(Y^{(n)},Z^{(n)}, K^{n} )$ is a  Cauchy sequences in $L^{2}_{\mF}(0, T+\delta; \mR)\times L^{2}_{\mF}(0, T; \mR^{d})\times H^{2}_{\mF}(0, T+\delta; \mR).$ Denote their limits by $(Y,Z, K).$  Taking limits in the above iterative equations, we obtain that $(Y,Z, K)$ satisfies the following RABSDE:
\begin{align*}
 \begin{cases}
  Y_t= \xi^{(2)}_T+\int^{T}_{t}f_{2}(s, Y_s, Z_s, Y_{s+\mu(s)},G_{s}(K), \mE^{\sF_{s}}[G_{s+\epsilon(s)}(K)])\dif s + K_{T}- K_{t}- \int^{T}_{t}Z_s\dif W_s,   t\in[0, T], \\
 Y_t= \xi^{(2)}_t, t \in [T,T+\delta]; K_t= \zeta^{(2)}_t,   t\in (T, T+\delta].\\
  \end{cases}
 \end{align*}
According to Theorem 3.1, we know
 $$Y_t = Y^{(2)}_t, K_t = K^{(2)}_t.$$
 Since  $Y^{(1)}_t \geq Y^{(3)}_t \geq Y^{(4)}_t \geq Y_t, K^{(1)}_t \leq K^{(3)}_t \leq K^{(4)}_t \leq K_t,\,\,a.e.,\,\,a.s.,$  it holds immediately that
 $$Y^{(1)}_t \geq Y^{(2)}_t, K^{(1)}_t \leq K^{(2)}_t   ,\,\,a.e.,\,\,a.s.$$

\end{proof}

\subsection{An application}In this section,  an application is given  for the above comparison theorem. We will prove that the following equation has a minimal solution with certain conditions.
\begin{align}
 \begin{cases}
 -\dif Y_{t}= f(t, Y_{t}, Z_{t}, \mE^{\sF_{t}}[Y_{t+\mu(t)}], G_{t}(K), \mE^{\sF_{t}}[G_{t+\epsilon(t)}(K)])\dif t +\dif K_{t}- Z_{t}\dif W_{t},  t\in[0, T], \\
  Y_{t}\geq S_{t}, t\in [0,T], \int^{T}_{0}(Y_{t}-S_{t})\dif K_{t}=0;\\
 Y_{t}= \xi_{t},t\in [T, T+\delta];  K_{t}=\zeta_{t},   t\in(T, T+\delta].
\end{cases}
\end{align}
We  make the following assumptions.
\begin{enumerate}[(B1)]
\item There exists a functional $\hat{f}$ such that $f(t,y,z,\theta, m, \bar{m})=\hat{f}(t,y,z,\mE^{\sF_{t}}[\theta], m, \bar{m})$ for any $t\in [0,T], y, m, \bar{m}\in \mR, z\in\mR^{d},\theta \in L^{2}(\sF_{r}; \mR), r\in [t, T+K].$

\item For any $t\in [0, T],y, m, \bar{m}\in \mR, z\in \mR^{d},$ $f(t,y,z,\cdot,m, \bar{m})$ is increasing, that is, $ f(t,y,z,\theta,m, \bar{m})\geq f(t,y,z,\theta^{\prime},m, \bar{m}),$ if $\theta \geq \theta^{\prime},$ $\theta, \theta^{\prime}\in L^{2}(\sF_{r}; \mR),\,\,r \in [t, T+K].$
\item  For any $t\in [0, T], \theta \in L^{2}(\sF_{r}; \mR), r \in [t, T+K] ,\bar{m},y\in \mR ,z\in \mR^{d},$ $f(t,y,z,\theta,\cdot, \bar{m})$ is decreasing, that is $f(t,y,z,\theta, m,\bar{m}) \leq f(t,y,z,\theta, m', \bar{m}),$ if  $m \geq m^{\prime}, m , m^{\prime}\in \mR.$

\item  For any $t\in [0, T], \theta \in L^{2}(\sF_{r}; \mR), r \in [t, T+K], m ,y\in \mR ,z\in \mR^{d},$ $f(t,y,z,\theta, m,\cdot)$ is decreasing, that is $f(t,y,z,\theta, m,\bar{m}) \leq f(t,y,z, \theta, m,\bar{ m}'),$ if  $\bar{m} \geq \bar{m}^{\prime}, \bar{m} , \bar{m}^{\prime}\in \mR.$
\item  For any $t\in [0, T], \theta \in L^{2}(\sF_{r}; \mR), r \in [t, T+K], y, m,\bar{m},m',\bar{m}'\in \mR ,z\in \mR^{d},$ we have
$$|f(t,y,z,\theta,m, \bar{m})|\leq C(h_{t}+|y|+|z|+\mE^{\sF_{t}}[|\theta|])+C_{1}(|m|+|\bar{m}|),$$
\begin{align*}
    |&f(t, y,z,\theta, m, \bar{m }) - f(t, y,z,\theta, m', \bar{m}')|
    \leq C_{1}(|m- m^{\prime}|+|\bar{m}- \bar{m}^{\prime}|).
    \end{align*}
\item For any  $t\in [0,T], r\in [t, T+\delta],$ $f(t,y,z,\theta,m, \bar{m})$ is continuous in $\mR\times\mR^{d}\times L^{2}(\sF_{r}; \mR)\times \mR\times \mR.$

\item If $y^{(1)}, y^{(2)}\in \sD(\mR)$ and $y^{(1)}_{t}\geq y^{(2)}_{t},$  we have $G_{t}(y^{(1)})\geq G_{t}(y^{(2)}).$
\end{enumerate}
\begin{enumerate}[(B8)]
\item If $y^{(n)}, y\in \sD(\mR)$ and $y^{(n)}_{t}\downarrow y_{t},$  we have $G_{t}(y^{(n)})\downarrow G_{t}(y).$
\end{enumerate}

\begin{enumerate}[(B9)]
\item For any $t \in [0, T], y,y^{\prime}\in \sD(\mR),$ $G_{t}$ satisfies $G(0)=0$  and
$$G_{t}(y)=G_{t}(\{y_{s\wedge t}\}_{0\leq s \leq T}), $$
$$\mE\bigg[\bigg. \int^{T}_{0}|G_{t}(y)-G_{t}(y')|^{2}\dif t \bigg]\bigg. \leq C_{1}\mE\bigg[\bigg. \int^{T}_{0}\bigg[\bigg.y_{t}-y'_{t}|^{2}\dif t\bigg]\bigg..$$
\end{enumerate}
\br
Obviously, though the condition $(\mathrm{B}9)$ is stronger than the condition$(\mathrm{H}4)$ , there also exist many examples satisfying $(\mathrm{B}9).$
$$ 1^{\circ} G_{t}: y\rightarrow \int^{t}_{\frac{t}{2}}y^{+}_{s}\dif s,2^{\circ} G_{t}: y\rightarrow y^{+}_{(t-\varepsilon)^{+}},3^{\circ}G_{t}: y\rightarrow y_{(t-\varepsilon)^{+}},4^{\circ} G_{t}: y\rightarrow \int^{t}_{\frac{t}{2}}y^{+}_{s}\dif s, 5^{\circ} G_{t}: y\rightarrow \frac{1}{\varepsilon}\int^{t}_{(t-\varepsilon)^{+}}y^{+}_{s}\dif s.$$
\er
Before giving our main results, we know that the following two equations have a unique triple of solution $(U^{(i)},V^{(i)}, \bar{K}^{(i)}),$ $i=1,2,$ respectively.
\begin{align*}
 \begin{cases}
 -\dif U^{(1)}_t= [C(h_{t}+| U^{(1)}_t|+ |V^{(1)}_t|+\mE^{\sF_{t}}[|U^{(1)}_{t+\mu(t)}|])+C_{1}(G_{t}(\bar{K}^{(1)})+ \mE^{\sF_{t}}[G_{t+\epsilon(t)}(\bar{K}^{(1)})])]\dif t +\dif \bar{K}^{(1)}_{t} - V^{(1)}_{t}\dif W_{t},\\
 \quad \quad \quad \quad \quad \quad \quad \quad\quad \quad \quad \quad \quad \quad \quad \quad\quad \quad \quad \quad \quad \quad \quad \quad \quad \quad \quad \quad \quad \quad \quad \quad \quad \quad \quad \quad \quad \quad \quad \quad t\in [0,T],\\
U^{(1)}_{t}\geq S_{t}, t\in [0,T];  \int^{T}_{0}(U^{(1)}_{t}-S_{t})\dif \bar{K}^{(1)}_{t}=0;\\
 U^{(1)}_{t}= \xi_{t}, t\in [T, T+\delta]; \bar{K}^{(1)}_{t}=\zeta_{t},  t\in (T, T+\delta],\\
 \end{cases}
\end{align*}
and
\begin{align*}
 \begin{cases}
 -\dif U^{(2)}_t= -[C(h_{t}+| U^{(2)}_t|+ |V^{(2)}_t|+\mE^{\sF_{t}}[|U^{(2)}_{t+\mu(t)}|])+C_{1}(G_{t}(\bar{K}^{(2)})+ \mE^{\sF_{t}}[G_{t+\epsilon(t)}(\bar{K}^{(2)})])]\dif t + \dif \bar{K}^{(2)}_{t}- V^{(2)}_{t}\dif W_{t},\\
 \quad \quad \quad \quad \quad \quad \quad \quad\quad \quad \quad \quad \quad \quad \quad \quad\quad \quad \quad \quad \quad \quad \quad \quad \quad \quad \quad \quad \quad \quad \quad \quad \quad \quad \quad \quad \quad \quad \quad \quad t\in [0,T],\\
U^{(2)}_{t}\geq S_{t}, t\in [0,T]; \int^{T}_{0}(U^{(2)}_{t}-S_{t})\dif \bar{K}^{(2)}_{t}=0;\\
 U^{(2)}_{t}= \xi_{t},t\in[T, T+\delta]; \bar{K}^{(2)}_{t}=\zeta_{t}, t\in(T, T+\delta].\\
 \end{cases}
\end{align*}
For any $t\in [0, T], r \in [t, T+K], \theta \in L^{2}(\sF_{r}; \mR),  y, m,\bar{m}\in \mR ,z\in \mR^{d},$ set
\begin{align*}
 f_{n}(t,  y, z,\theta, m, \bar{m})&=\inf_{(a,b,\vartheta)\in \mR\times \mR^{d}\times  L^{2}(\sF_{r}; \mR)}\big\{\big. f(t,  a, b,\vartheta,  m, \bar{m})+n|y-a|+ n|z-b| +n\mE^{\sF_{t}}[|\theta-\vartheta| ]\big\}\big..
\end{align*}
 The following lemma is mainly from Lemma 2.8 and Lemma 3.5 in Yang and Elliot (2013).

\bl
Assume that $(\mathrm{B}1)-(\mathrm{B}9)$ are established. we have the following properties.
\begin{enumerate}[($\mathrm{i})$]
\item  Linear growth: for any $t\in [0, T],  y,m,\bar{m}\in \mR, z\in \mR^{ d},  r\in [t, T+K], \theta\in L^{2}(\sF_{r}; \mR),$ we have
$$|f_{n}(t,  y, z, \theta ,m,\bar{m} )|\leq C(h_{t}+|y|+ |z|+\mE^{\sF_{t}}[|\theta|])+ C'(|m|+|\bar{m}|). $$
\end{enumerate}
\begin{enumerate}[($\mathrm{ii})$]
\item Monotonicity in $n$: for any $t\in [0, T], y\in \mR, z\in \mR^{ d}, m, \bar{m}\in \mR ,r\in [t, T+\delta], \theta\in L^{2}(\sF_{r}; \mR),$ we have
$f_{n}(t,  y, z, \theta, m,\bar{m})$ is increasing in $n$.
\end{enumerate}

\begin{enumerate}[($\mathrm{iii})$]
\item Monotonicity in $\theta:$ for any $t\in [0, T],  y\in \mR, m, \bar{m}\in \mR ,z\in \mR^{ d}$, we know $f_{n}(t,  y, z, \theta, m,\bar{m} )$ is increasing in $\theta,$ that is $f_{n}(t,  y, z,\theta_{1}, m,\bar{m})\geq f_{n}(t,  y, z, \theta_{2},m,\bar{m})$ if $\theta_{1} \geq \theta_{2}, r\in [t, T+K], \theta_{1},\theta_{2}\in L^{2}(\sF_{r}; \mR).$

\end{enumerate}

\begin{enumerate}[($\mathrm{iv})$]

\item  For any $t\in [0, T], \theta \in L^{2}(\sF_{r}; \mR), r \in [t, T+K] ,\bar{m},y\in \mR ,z\in \mR^{d},$ $f(t,y,z,\theta,\cdot, \bar{m})$ is decreasing, that is $f_{n}(t,y,z,\theta, m,\bar{m}) \leq f_{n}(t,y,z,m', \theta, \bar{m}),$ if  $m \geq m^{\prime}, m , m^{\prime}\in \mR.$

\end{enumerate}

\begin{enumerate}[($\mathrm{v})$]
\item  For any $t\in [0, T], \theta \in L^{2}(\sF_{r}; \mR), r \in [t, T+K], m ,y\in \mR ,z\in \mR^{d},$ $f(t,y,z,\theta, m,\cdot)$ is decreasing, that is $f_{n}(t,y,z,\theta, m,\bar{m}) \leq f_{n}(t,y,z, \theta, m,\bar{ m}'),$ if  $\bar{m} \geq \bar{m}^{\prime}, \bar{m} , \bar{m}^{\prime}\in \mR.$

\end{enumerate}

\begin{enumerate}[($\mathrm{vi})$]
\item Lipschitz condition: for any $t\in [0, T],  y,y^{\prime},m,m^{\prime},\bar{m},\bar{m}^{\prime}\in \mR, z\in \mR^{ d}, r\in [t, T+\delta], \theta,  \theta^{\prime}\in L^{2}(\sF_{r}; \mR), $ we have
    \begin{align*}
    |f_{n}(t,  y, z,\theta, m,\bar{m}) - f_{n}(t,  y', z',\theta^{\prime}, m',\bar{m}')|
    \leq n(|y- y^{\prime}| + |z-z^{\prime}| + \mE^{\sF_{t}}[|\theta - \theta^{\prime}|])+  C_{1}|m- m^{\prime}|+C_{1}|\bar{m}- \bar{m}^{\prime}|.
     \end{align*}
\end{enumerate}
\begin{enumerate}[($\mathrm{vii})$]
\item Convergence: for any   $t\in [0, T], r\in [t, T+K],$ if $(y_{n}, z_{n},\theta^{(n)}, m_{n},\bar{m}_{n} )\rightarrow (y,z, \theta,m,\bar{m})\,\, \mbox{in}\,\, \mR\times\mR^{d}\times L^{2}(\sF_{r}; \mR)\times\mR\times\mR, n\rightarrow \infty,$ then there exists a subsequence if necessary such that
$$f_{n}(t,y_{n},z_{n}, \theta^{(n)}, m_{n},\bar{m}_{n})\rightarrow  f(t,u,y,z,\theta,m, \bar{m}), n\rightarrow \infty.$$
\end{enumerate}
\el
For each $n,$ the following equation has a unique triple of solution $(Y^{(n)}, Z^{(n)}, K^{(n)}).$
\begin{align}
 \begin{cases}
 -\dif Y^{(n)}_{t}= f_{n}(t, Y^{(n)}_{t}, Z^{(n)}_{t}, Y^{(n)}_{t+\mu(t)}, G^{(n)}_{t}(K), \mE^{\sF_{t}}[G_{t+\epsilon(t)}(K^{(n)})])\dif t+\dif K^{(n)}_{t} - Z^{(n)}_{t}\dif W_{t},  t\in[0, T], \\
  Y^{(n)}_{t}\geq S_{t}, t\in [0,T], \int^{T}_{0}(Y^{(n)}_{t}-S_{t})\dif K^{(n)}_{t}=0;\\
 Y^{(n)}_{t}= \xi_{t},t\in [T, T+\delta];  K^{(n)}_{t}=\zeta_{t},   t\in(T, T+\delta].
\end{cases}
\end{align}
\bl
Set $\lambda=24C^{2},\beta=\lambda+\frac{6C^{2}(1+L)}{\lambda}+2, \gamma=4[((4+4e^{\beta T})12Te^{\beta T}+4Te^{\beta T})6C^{2}(1+L)+4e^{\beta T}]$. If  the following conditions
\begin{enumerate}[($\mathrm{a})$]
\item $C_{1}$ is small enough such that $((4+4e^{\beta T})12Te^{\beta T}+4Te^{\beta T})6C^{2}_{1}(T+L(T+\delta))+\frac{1+6C_{1}^{2}(T+L(T+\delta))\gamma}{\lambda}< 1.$
\end{enumerate}

\begin{enumerate}[($\mathrm{b})$]
\item  $(\mathrm{B}1),(\mathrm{B}5),(\mathrm{B}9)$ are established.
\end{enumerate}

\begin{enumerate}[($\mathrm{c})$]
\item  $\xi, \zeta\in S^{2}_{\mF}(T, T + \delta; \mR) $ and $S \in S^{2}_{\mF}(0, T; \mR) .$

\end{enumerate}
 hold, then there exists a constant $M>0$ such that
$$\sup_{n}\mE\bigg[\bigg.\sup_{t\in [0, T+\delta]}|Y^{(n)}_{t}|^{2}+\sup_{t\in [0, T+\delta]}|K^{(n)}_{t}|^{2}+\int^{T}_{0}|Z^{(n)}_{t}|^{2}\dif t\bigg]\bigg. \leq M.$$
\el
\begin{proof}

By It\^{o}'s formula, we have
\begin{align*}
 &\mE[e^{\beta t}|Y^{(n)}_t|^{2}] +\beta\mE\bigg[\bigg.\int^{T}_{t}e^{\beta s}|Y^{(n)}_{s}|^{2}\dif s\bigg]\bigg.+\mE\bigg[\bigg.\int^{T}_{t}e^{\beta s}|Z^{(n)}_{s}|^{2}\dif s\bigg]\bigg. \\
 &=\mE[e^{\beta T}|\xi_T|^{2}]+  2\mE\bigg[\bigg.\int^{T}_{t}e^{\beta s} Y^{(n)}_s f_{n}(s,Y^{(n)}_s, Z^{(n)}_s, Y^{(n)}_{s+\mu(s)},G_s(K^{(n)}), \mE^{\sF_{s}}[G_{s+\epsilon(s)}(K^{(n)})] ) \dif s\bigg]\bigg. +2\mE\bigg[\bigg.\int^{T}_{t}e^{\beta s} Y^{(n)}_s  \dif K^{(n)}_s  \bigg]\bigg..
\end{align*}
Using Young's inequality, we have
\begin{align*}
  &\mE[e^{\beta t}|Y^{(n)}_t|^{2}] +\beta\mE\bigg[\bigg.\int^{T}_{t}e^{\beta s}|Y^{(n)}_{s}|^{2}\dif s\bigg]\bigg.+\mE\bigg[\bigg.\int^{T}_{t}e^{\beta s}|Z^{(n)}_{s}|^{2}\dif s\bigg]\bigg. \\
 &=\mE[e^{\beta T}|\xi_T|^{2}]+  2\mE\bigg[\bigg.\int^{T}_{t}e^{\beta s} Y^{(n)}_s f_{n}(s,Y^{(n)}_s, Z^{(n)}_s, Y^{(n)}_{s+\mu(s)},G_s(K^{(n)}), \mE^{\sF_{s}}[G_{s+\epsilon(s)}(K^{(n)})] ) \dif s\bigg]\bigg. +2\mE\bigg[\bigg.\int^{T}_{t}e^{\beta s} Y^{(n)}_s  \dif K^{(n)}_s  \bigg]\bigg.\\
 &\leq \mE[e^{\beta T}|\xi_T|^{2}]+  2\mE\bigg[\bigg.\int^{T}_{t} e^{\beta s}|Y^{(n)}_s| |f^{(n)}(s, Y^{(n)}_s, Z^{(n)}_s, Y^{(n)}_{s+\mu(s)},G_s(K^{(n)}), \mE^{\sF_{s}}[G_{s+\epsilon(s)}(K^{(n)})])| \dif s\bigg]\bigg.+2\mE\bigg[\bigg.\int^{T}_{t}e^{\beta s} Y^{(n)}_s  \dif K^{(n)}_s  \bigg]\bigg.\\
 & \leq \mE[e^{\beta T}|\xi_T|^{2}]+  2\mE\bigg[\bigg.\int^{T}_{t}e^{\beta s} |Y^{(n)}_s|(Ch_{s} +C|Y^{(n)}_s|+C|Z^{(n)}_{s}|+C|Y^{(n)}_{s+\mu(s)}|\\
 & \quad \quad \quad \quad\quad \quad \quad \quad\quad \quad \quad \quad \quad\quad \quad \quad \quad + C_{1}|G_s(K^{(n)})|+C_{1}\mE^{\sF_{s}}[|G_{s+\epsilon(s)}(K^{(n)})|])\dif s\bigg]\bigg.+2\mE\bigg[\bigg.\int^{T}_{t}e^{\beta s} Y^{(n)}_s  \dif K^{(n)}_s  \bigg]\bigg.\\
 &\leq C^{\prime}+\bigg(\bigg.\lambda +\frac{6C^{2}(1+L)}{\lambda}  \bigg)\bigg.\mE\bigg[\bigg.\int^{T}_{t}e^{\beta s} |Y^{(n)}_s|^{2}\dif s\bigg]\bigg.+\frac{6C^{2}}{\lambda}\mE\bigg[\bigg.\int^{T}_{t}e^{\beta s} |Z^{(n)}_s|^{2}\dif s\bigg]\bigg.+\frac{6C^{2}_{1}(T+L(T+\delta))}{\lambda}\mE\bigg[\bigg. \sup_{0\leq t\leq T}e^{\beta s} |K^{(n)}_s|^{2}\dif s\bigg]\bigg. \\
 & \quad + 2\mE\bigg[\bigg.\int^{T}_{t} e^{\beta s}Y^{(n)}_s  \dif K^{(n)}_s  \bigg]\bigg..
\end{align*}
Since,
$$2\mE\bigg[\bigg.\int^{T}_{0}e^{\beta s} Y^{(n)}_s  \dif K^{(n)}_s  \bigg]\bigg. \leq \lambda\gamma\mE\bigg[\bigg.\sup_{0\leq t\leq T} e^{\beta t} |S_t|^{2}\dif s\bigg]\bigg.+ \frac{1}{\lambda\gamma}\mE\bigg[\bigg. \sup_{0\leq t\leq T}e^{\beta t} |K^{(n)}_s|^{2}\dif s\bigg]\bigg.,$$
It leads to
\begin{align}
 &\mE[e^{\beta t}|Y^{(n)}_t|^{2}] +\beta\mE\bigg[\bigg.\int^{T}_{t}e^{\beta s}|Y^{(n)}_{s}|^{2}\dif s+\mE\bigg[\bigg.\int^{T}_{t}e^{\beta s}|Z^{(n)}_{s}|^{2}\dif s\bigg]\bigg.\no\\
 &\leq C^{\prime}+\bigg(\bigg.\lambda +\frac{6C^{2}(1+L)}{\lambda}  \bigg)\bigg.\mE\bigg[\bigg.\int^{T}_{t}e^{\beta s} |Y^{(n)}_s|^{2}\dif s\bigg]\bigg.+\frac{6C^{2}}{\lambda}\mE\bigg[\bigg.\int^{T}_{t}e^{\beta s} |Z^{(n)}_s|^{2}\dif s\bigg]\bigg.\no\\
 & \quad +\bigg(\bigg.\frac{6C^{2}_{1}(T+L(T+\delta))}{\lambda}+\frac{1}{\lambda\gamma}\bigg)\bigg.\mE\bigg[\bigg. \sup_{0\leq t\leq T}e^{\beta s} |K^{(n)}_s|^{2}\dif s\bigg]\bigg..
  \end{align}
Moreover,
\begin{align*}
K^{(n)}_t=-Y^{(n)}_{t}+Y^{(n)}_{0}-\int^{t}_{0}f_{n}(s,Y^{(n)}_s, Z^{(n)}_s, Y^{(n)}_{s+\mu(s)},G_s(K^{(n)}), \mE^{\sF_{s}}[G_{s+\epsilon(s)}(K^{(n)})]) \dif s + \int^{t}_{0}Z^{(n)}_{s}\dif W_{s}.
\end{align*}
Then,

\begin{align*}
\frac{1}{\gamma}e^{\beta t}|K^{(n)}_{t}|^{2}&\leq e^{\beta t}\frac{4}{\gamma}|Y^{(n)}_{t}|^{2}+\frac{4}{\gamma}|e^{\beta t}Y^{(n)}_{0}|^{2}+\frac{4T}{\gamma}e^{\beta t}\int^{t}_{0}|f_{n}(s,Y^{(n)}_s, Z^{(n)}_s, Y^{(n)}_{s+\mu(s)},G_s(K^{(n)}),, \mE^{\sF_{s}}[G_{s+\epsilon(s)}(K^{(n)})]) |^{2}\dif s \\
&\quad + \frac{4}{\gamma}e^{\beta t}\bigg|\bigg. \int^{t}_{0}\hat{Z}^{(n)}_{s}\dif W_{s}\bigg|\bigg.^{2}.
\end{align*}
By BDG's formula, we have
\begin{align}
\mE[\sup_{0\leq t\leq T}\frac{1}{\gamma}e^{\beta t}|K^{(n)}_{t}|^{2}]&\leq\frac{4+4e^{\beta T}}{\gamma}\mE[\sup_{0\leq t\leq T}e^{\beta t}|Y^{(n)}_{t}|^{2}] + \frac{4e^{\beta T}}{\gamma}\mE\bigg[\bigg. \int^{T}_{0}e^{\beta s}|Z^{(n)}_{s}|^{2}\dif s\bigg]\bigg. \no\\
&\quad+\frac{4 Te^{\beta T}}{\gamma}\int^{T}_{0}e^{\beta s}|f_{n}(s,Y^{(n)}_s, Z^{(n)}_s, Y^{(n)}_{s+\mu(s)},G_s(K^{(n)}), \mE^{\sF_{s}}[G_{s+\epsilon(s)}(K^{(n)})]) |^{2}\dif s
\end{align}
With the following representation,
\begin{align*}
Y^{(n)}_{t}=\mathrm{ess}\sup_{\tau \in \Gamma_{t}}\mE\bigg[\bigg. \int^{\tau}_{t} f_{n}(s,  Y^{(n)}_s, Z^{(n)}_s, Y^{(n)}_{s+\mu(s)}, G_{s}(K^{(n)}), \mE^{\sF_{s}}[G_{s+\epsilon(s)}(K^{(n)})])\dif s + S_{\tau}1_{\tau<T} +\xi_{T}1_{\tau=T}\bigg| \bigg. \sF_{t}   \bigg]\bigg.,
\end{align*}
we have
\begin{align*}
|Y^{(n)}_{t}|\leq \mE\bigg[\bigg. \int^{T}_{0} |f_{n}(s,  Y^{(n)}_s, Z^{(n)}_s, Y^{(n)}_{s+\mu(s)}, G_{s}(K^{(n)}) , \mE^{\sF_{s}}[G_{s+\epsilon(s)}(K^{(n)})])|\dif s+ \sup_{0\leq t\leq T}|S_{t}| +|\xi_{T}| \bigg| \bigg. \sF_{t}   \bigg]\bigg..
\end{align*}
It leads to
\begin{align*}
|Y^{(n)}_{t}|^{2}
 \leq 3\bigg(\bigg. \mE\bigg[\bigg. \int^{T}_{0} |f_{n}(s,  Y^{(n)}_s, Z^{(n)}_s, Y^{(n)}_{s+\mu(s)}, G_{s}(K^{(n)}), \mE^{\sF_{s}}[G_{s+\epsilon(s)}(K^{(n)})])|\dif s \bigg| \bigg. \sF_{t}   \bigg]\bigg.  \bigg)\bigg.^{2}+3\bigg(\bigg. \mE\bigg[\bigg.  \sup_{0\leq t\leq T}|S_{t}| \bigg| \bigg. \sF_{t}   \bigg]\bigg.  \bigg)\bigg.^{2}+3\bigg(\bigg. \mE[ |\xi_{T}||  \sF_{t}   ]  \bigg)\bigg.^{2}.
\end{align*}
Doob's maximal inequality implies that
\begin{align*}
\mE\bigg[\bigg.\sup_{0\leq t \leq T}|Y^{(n)}_t|^{2}\bigg]\bigg. \leq 12T \mE\bigg[\bigg. \int^{T}_{0} |f_{n}(s,  Y^{(n)}_s, Z^{(n)}_s, Y^{(n)}_{s+\mu(s)}, G_{s}(K^{(n)}), \mE^{\sF_{s}}[G_{s+\epsilon(s)}(K^{(n)})] )|^{2}\dif s    \bigg]\bigg.  +12 \mE\bigg[\bigg.  \sup_{0\leq t\leq T}|S_{t}|^{2}    \bigg]\bigg.  +12 \mE[ |\xi_{T}|^{2} ].
 \end{align*}
Then, we have

 \begin{align}
\mE\bigg[\bigg.\sup_{0\leq t \leq T}\frac{1}{\gamma}e^{\beta t}|Y^{(n)}_t|^{2}\bigg]\bigg. &\leq  \frac{12Te^{\beta T}}{\gamma}\mE\bigg[\bigg.\int^{T}_{0} e^{\beta s}|f_{n}(s,  Y^{(n)}_s, Z^{(n)}_s, Y^{(n)}_{s+\mu(s)}, G_{s}(K^{(n)}), \mE^{\sF_{s}}[G_{s+\epsilon(s)}(K^{(n)})])|^{2}\dif s \bigg]\bigg. + C'.
\end{align}
Combining (9) and  (10), we derive
\begin{align}
\mE[\sup_{0\leq t\leq T}\frac{1}{\gamma}e^{\beta t}|K^{(n)}_{t}|^{2}]&\leq\frac{(4+4e^{\beta T})12Te^{\beta T}+4Te^{\beta T}}{\gamma} \mE\bigg[\bigg.\int^{T}_{0} e^{\beta s}|f_{n}(s,  Y^{(n)}_s, Z^{(n)}_s, Y^{(n)}_{s+\mu(s)}, G_{s}(K^{(n)}), \mE^{\sF_{s}}[G_{s+\epsilon(s)}(K^{(n)})])|^{2}\dif s\bigg]\bigg. \no\\
&\quad + \frac{4e^{\beta T}}{\gamma}\mE\bigg[\bigg. \int^{T}_{0}e^{\beta s}|Z^{(n)}_{s}|^{2}\dif s\bigg]\bigg.+C'.
\end{align}
From (8) and  (11), we obtain
\begin{align*}
 \mE&\bigg[\bigg. \sup_{0\leq t\leq T} \frac{1}{\gamma}e^{\beta s}|K^{(n)}_s|^{2}\dif s\bigg]\bigg. +\beta\mE\bigg[\bigg.\int^{T}_{t}e^{\beta s}|Y^{(n)}_{s}|^{2}\dif s\bigg]\bigg. +\mE\bigg[\bigg.\int^{T}_{t}e^{\beta s}|Z^{(n)}_{s}|^{2}\dif s\bigg]\bigg.\\
 &\leq \bigg(\bigg.\lambda +\frac{6C^{2}(1+L)}{\lambda}  \bigg)\bigg.\mE\bigg[\bigg.\int^{T}_{t}e^{\beta s} |Y^{(n)}_s|^{2}\dif s\bigg]\bigg.+\bigg(\bigg.\frac{6C^{2}}{\lambda}+\frac{4e^{\beta T}}{\gamma}\bigg)\bigg. \mE\bigg[\bigg.\int^{T}_{t}e^{\beta s} |Z^{(n)}_s|^{2}\dif s\bigg]\bigg.\\
&\quad  +\frac{(4+4e^{\beta Y})12Te^{\beta T}+4Te^{\beta T}}{\gamma} \mE\bigg[\bigg.\int^{T}_{0} e^{\beta s}|f_{n}(s,  Y^{(n)}_s, Z^{(n)}_s, Y^{(n)}_{s+\mu(s)}, G_{s}(K^{(n)}),\mE^{\sF_{s}}[G_{s+\epsilon(s)}(K^{(n)})])|^{2}\dif s\bigg]\bigg. \no\\
&\quad + \bigg(\bigg.\frac{1}{\lambda\gamma} +\frac{6C_{1}^{2}(T+L(T+\delta))}{\lambda}  \bigg)\bigg.\mE\bigg[\bigg.  \sup_{0\leq t\leq T}e^{\beta t}|K^{(n)}_{t}|^{2}     \bigg]\bigg.+C'.
\end{align*}
From $(\mathrm{i})$ in Lemma 3.5, we have
\begin{align*}
 &\mE\bigg[\bigg. \sup_{0\leq t\leq T} \frac{1}{\gamma}e^{\beta s}|K^{(n)}_s|^{2}\dif s\bigg]\bigg. +\beta\mE\bigg[\bigg.\int^{T}_{t}e^{\beta s}|Y^{(n)}_{s}|^{2}\dif s+\mE\bigg[\bigg.\int^{T}_{t}e^{\beta s}|Z^{(n)}_{s}|^{2}\dif s\bigg]\bigg.\\
& \leq \bigg(\bigg.\lambda +\frac{6C^{2}(1+L)}{\lambda}+\frac{((4+4e^{\beta T})12Te^{\beta T}+4Te^{\beta T})6C^{2}(1+L)}{\gamma}   \bigg)\bigg.\mE\bigg[\bigg.\int^{T}_{t}e^{\beta s} |Y^{(n)}_s|^{2}\dif s\bigg]\bigg.\\
&\quad +\bigg(\bigg.\frac{((4+4e^{\beta T})12Te^{\beta T}+4Te^{\beta T})6C^{2}}{\gamma}+ \frac{6C^{2}}{\lambda}+\frac{4e^{\beta T}}{\gamma}\bigg)\bigg.\mE\bigg[\bigg.\int^{T}_{t}e^{\beta s} |Z^{(n)}_s|^{2}\dif s\bigg]\bigg.\\
&\quad+ \bigg(\bigg.((4+4e^{\beta T})12Te^{\beta T}+4Te^{\beta T})6C^{2}_{1}(T+L(T+\delta))+\frac{1+6C_{1}^{2}(T+L(T+\delta))\gamma}{\lambda}\bigg)\bigg.\mE\bigg[\bigg. \sup_{0\leq t\leq T}\frac{1}{\gamma}e^{\beta t}|K^{(n)}_{t}|^{2}    \bigg]\bigg.+ C^{\prime}.
 \end{align*}
Choosing $\lambda=24C^{2},\beta=\lambda+\frac{6C^{2}(1+L)}{\lambda}+2, \gamma=4[((4+4e^{\beta T})12Te^{\beta T}+4Te^{\beta T})6C^{2}(1+L)+4e^{\beta T}],$   then there exists a constant $M'>0$  such that
\begin{align*}
 \mE\bigg[\bigg. \sup_{0\leq t\leq T} \frac{1}{\gamma}e^{\beta s}|K^{(n)}_s|^{2}\dif s\bigg]\bigg. +\mE\bigg[\bigg.\int^{T}_{0}e^{\beta s}|Y^{(n)}_{s}|^{2}\dif s +\mE\bigg[\bigg.\int^{T}_{0}e^{\beta s}|Z^{(n)}_{s}|^{2}\dif s \bigg]\bigg. \leq M^{\prime}.
\end{align*}
The proof is complete.
\end{proof}

\bt
 For Eq.(6),  let the assumptions  $(\mathrm{B}1)-(\mathrm{B}9), (\mathrm{i}),(\mathrm{ii})$  be in force.   if $\xi \in S^{2}_{\mF}(T,T+K; \mR),\zeta \in H^{2}_{\mF}(T,T+K; \mR),S \in S^{2}_{\mF}(0,T; \mR)$ and $C_{1}> 0$  is small enough,  then there exists a minimal solution $Y$, that is if $\tilde{Y}$ is another solution of Eq.(6),  we have
 $$ Y_{t}\leq \tilde{Y}_{t},a.e.,a.s. $$
 \et

\begin{proof}
By  Theorem 3.3,  for any $m<n,$ we have $U^{(2)}_{t}\leq Y^{(m)}_{t}\leq Y^{(n)}_{t}\leq U^{(1)}_{t}, \bar{K}^{(1)}_{t}\leq K^{(n)}_{t}\leq K^{(m)}_{t}\leq \bar{K}^{(2)}_{t}, a.e., a.s.$ Thus, there exist a process $\{Y_{t}, t\in [0, T+\delta]\}$  and  a process $\{K_{t}, t\in [0, T+\delta]\}$  such that $Y^{(n)}_{t}\uparrow Y_{t}, K^{(n)}_{t}\downarrow K_{t},  n\rightarrow \infty.$ From monotone convergence theorem, we have $\lim_{n\rightarrow \infty}\mE\int^{T+\delta}_{0}|Y^{(n)}_{t}-Y_{t}|^{2}\dif t=0,$ $\lim_{n\rightarrow \infty}\mE\int^{T+\delta}_{0}|K^{(n)}_{t}-K_{t}|^{2}\dif t=0.$ Using It\^{o}'s
formula to $|Y^{(n)}_{t}-Y^{(m)}_{t}|^{2},$ we obtain
\begin{align*}
 &\mE[|Y^{(n)}_{t}-Y^{(m)}_{t}|^{2}]+\mE\bigg[\bigg.\int^{T}_{t}|Z^{(n)}_{s}-Z^{(m)}_{s}|^{2}\dif s \bigg]\bigg.
 =  2\mE\bigg[\bigg.\int^{T}_{t} (Y^{(n)}_{s}-Y^{(m)}_{s})(f_{n}(s,  Y^{(n)}_s, Z^{(n)}_s, Y^{(n)}_{s+\mu(s)}, G_{s}(K^{(n)}), \mE^{\sF_{s}}[G_{s+\epsilon(s)}(K^{(n)})])\\
&\quad \quad \quad \quad \quad \quad \quad \quad \quad\quad \quad \quad \quad \quad \quad \quad\quad \quad -f_{m}(s, Y^{(m)}_s, Z^{(m)}_s, Y^{(m)}_{s+\mu(s)}), G_{s}(K^{(m)}),\mE^{\sF_{s}}[G_{s+\epsilon(s)}(K^{(m)})]) \dif s\bigg]\bigg..
 \end{align*}
Set
$$f_{m,n}(s)=f_{n}(s, Y^{(n)}_s,Z^{(n)}_s, Y^{(n)}_{s+\mu(s)}, G_{s}(K^{(n)}),\mE^{\sF_{s}}[G_{s+\epsilon(s)}(K^{(n)})])
-f_{m}(s, Y^{(m)}_s, Z^{(m)}_s, Y^{(m)}_{s+\mu(s)}, G_{s}(K^{(m)}), \mE^{\sF_{s}}[G_{s+\epsilon(s)}(K^{(m)})]).$$
From Lemma 3.6,  we have
\begin{align*}
&\mE\bigg[\bigg.\int^{T}_{t}|Z^{(n)}_{s}-Z^{(m)}_{s}|^{2}\dif s\bigg]\bigg.
\leq 2\bigg(\bigg. \mE\bigg[\bigg.\int^{T}_{t}|f_{m,n}(s)|^{2} \dif s\bigg]\bigg.\bigg)\bigg.^{\frac{1}{2}}
 \bigg(\bigg. \mE\bigg[\bigg.\int^{T}_{t}|Y^{(n)}_{s}-Y^{(m)}_{s}|^{2} \dif s\bigg]\bigg.\bigg)\bigg.^{\frac{1}{2}}\\
& \leq C'\bigg(\bigg.\mE\bigg[\bigg.\int^{T}_{t}|Y^{(n)}_{s}-Y^{(m)}_{s}|^{2} \dif s\bigg]\bigg.\bigg)\bigg.^{\frac{1}{2}}.
\end{align*}
Thus,
$$\lim_{n\rightarrow \infty}\mE\bigg[\bigg.\int^{T}_{0}|Z^{(n)}_{s}-Z^{(m)}_{s}|^{2}\dif s\bigg]\bigg. =0.$$
By Cauchy convergence criterion, there exists a process $\{Z_{s}, t\in [0, T]\}$ such that
 $$\lim_{n\rightarrow \infty}\mE\bigg[\bigg.\int^{T}_{0}|Z^{(n)}_{s}-Z_{s}|^{2}\dif s\bigg]\bigg. =0.$$
By It\^{o}'s formula,
\begin{align*}
 &|Y^{(n)}_{t}-Y^{(m)}_{t}|^{2}+\int^{T}_{t}|Z^{(n)}_{s}-Z^{(m)}_{s}|^{2}\dif s
 =  2\int^{T}_{t} (Y^{(n)}_{s}-Y^{(m)}_{s})f_{m,n}(s) \dif s
   + 2\int^{T}_{t}(Y^{(n)}_{s}-Y^{(m)}_{s})(Z^{(n)}_{s}-Z^{(m)}_{s})\dif W_{s}\\
 & \leq 2\bigg(\bigg. \int^{T}_{t}|f_{m,n}(s)|^{2} \dif s\bigg)\bigg.^{\frac{1}{2}} \bigg(\bigg.\int^{T}_{t}|Y^{(n)}_{s}-Y^{(m)}_{s}|^{2} \dif s\bigg)\bigg.^{\frac{1}{2}}\\
 &\quad  + 2\int^{T}_{t}(Y^{(n)}_{s}-Y^{(m)}_{s})(Z^{(n)}_{s}-Z^{(m)}_{s})\dif W_{s}.
 \end{align*}
By BDG's inequality, we have
\begin{align*}
\mE[\sup_{0\leq t\leq T}|Y^{(n)}_{t}-Y^{(m)}_{t}|^{2}]\leq C'\bigg(\bigg.\mE\bigg[\bigg.\int^{T}_{0}|Y^{(n)}_{s}-Y^{(m)}_{s}|^{2} \dif s\bigg.\bigg)\bigg.^{\frac{1}{2}}+ C'\mE\bigg[\bigg.\int^{T}_{0}|Z^{(n)}_{s}-Z^{(m)}_{s}|^{2}\dif s\bigg]\bigg.\rightarrow 0, n\rightarrow 0.
 \end{align*}
Thus,  $\mE[\sup_{0\leq t\leq T}|Y^{(n)}_{t}-Y_{t}|^{2}]\rightarrow 0, n\rightarrow 0.$
Since,
\begin{align*}
K^{(n)}_{t}-K^{(m)}_{t}&=-(Y^{(n)}_{t}-Y^{(m)}_{t})+Y^{(n)}_{0}-Y^{(m)}_{0}-\int^{t}_{0}f_{m,n}(s)\dif s
+ \int^{t}_{0}(Z^{(n)}_{s}-Z^{(m)}_{s})\dif W_{s},
\end{align*}
Thus,
\begin{align*}
|K^{(n)}_{t}-K^{(m)}_{s}|^{2}&\leq C'\int^{t}_{0}|f_{m,n}(s)|^{2}\dif s +C'\bigg|\bigg. \int^{t}_{0}(Z^{(n)}_{s}-Z^{(m)}_{s})\dif W_{s}\bigg|\bigg.^{2}
+C'|Y^{(n)}_{t}-Y^{(m)}_{t}|^{2}+C'|Y^{(n)}_{0}-Y^{(m)}_{0}|^{2}.
\end{align*}
Using BDG's inequality, we derive

\begin{align}
\mE\bigg[\bigg.\sup_{0\leq s\leq T}|K^{(n)}_{s}-K^{(m)}_{s}|^{2}\bigg]\bigg.&\leq  C'\mE\bigg[\bigg.\int^{T}_{0}|f_{m,n}(s)|^{2}\dif s\bigg]\bigg.
 +C'\mE\bigg[\bigg.\sup_{0\leq t\leq T}|Y^{(n)}_{s}-Y^{(m)}_{s}|^{2}\bigg]\bigg. + C'\mE\bigg[\bigg. \int^{T}_{0}|Z^{(n)}_{s}-Z^{(m)}_{s}|^{2}\dif s\bigg]\bigg..
\end{align}
Since,
\begin{align*}
Y^{(n)}_{t}=\mathrm{ess}\sup_{\tau \in \Gamma_{t}}\mE\bigg[\bigg. \int^{\tau}_{t} f_{n}(s,  Y^{(n)}_s, Z^{(n)}_s, Y^{(n)}_{s+\mu(s)}, G_{s}(K^{(n)}), \mE^{\sF_{s}}[G_{s+\epsilon(s)}(K^{(n)})])\dif s + S_{\tau}1_{\tau<T} +\xi_{T}1_{\tau=T}\bigg| \bigg. \sF_{t}   \bigg]\bigg.,
\end{align*}
and
\begin{align*}
Y^{(m)}_{t}=\mathrm{ess}\sup_{\tau \in \Gamma_{t}}\mE\bigg[\bigg. \int^{\tau}_{t} f_{m}(s,  Y^{(m)}_s, Z^{(m)}_s, Y^{(m)}_{s+\mu(s)}, G_{s}(K^{(m)}), \mE^{\sF_{s}}[G_{s+\epsilon(s)}(K^{(m)})])\dif s + S_{\tau}1_{\tau<T} +\xi_{T}1_{\tau=T}\bigg| \bigg. \sF_{t}   \bigg]\bigg.,
\end{align*}
using similar method, we derive
\begin{align}
\mE\bigg[\bigg.\sup_{0\leq t \leq T}|Y^{n}_{t}-Y^{m}_{t}|^{2}\bigg]\bigg. \leq  C'\mE\bigg[\bigg.\int^{T}_{0} |f_{m,n}(s)|^{2}\dif s\bigg]\bigg..
 \end{align}
(12) and (13) lead to
\begin{align}
&\mE\bigg[\bigg.\sup_{0\leq t\leq T}|K^{(n)}_{t}-K^{(m)}|^{2}\bigg]\bigg.\leq  C'\mE\bigg[\bigg.\int^{T}_{0}|f_{m,n}(s)|^{2}\dif s\bigg]\bigg.
  + C'\mE\bigg[\bigg. \int^{T}_{0}|Z^{(n)}_{t}-Z^{(m)}_{t}|^{2}\dif s\bigg]\bigg.\no\\
 &\leq    C'\mE\bigg[\bigg.\int^{T}_{0}|f_{n}(s, Y^{(n)}_s, Z^{(n)}_s, Y^{(n)}_{s+\mu(s)},G_{s}(K^{(n)}), \mE^{\sF_{s}}[G_{s+\epsilon(s)}(K^{(n)})])-f(s, Y_s, Z_s, Y_{s+\mu(s)},G_{s}(K), \mE^{\sF_{s}}[G_{s+\epsilon(s)}(K)])|^{2}\dif s\bigg]\bigg.\no\\
 &\quad + C'\mE\bigg[\bigg. \int^{T}_{0}|Z^{(n)}_{s}-Z^{(m)}_{s}|^{2}\dif s\bigg]\bigg.\no \\
  &\quad +C'\mE\bigg[\bigg.\int^{T}_{0}|f_{m}(s, Y^{(m)}_s, Z^{(m)}_s, Y^{(m)}_{s+\mu(s)},G_{s}(K^{(m)}) ,\mE^{\sF_{s}}[G_{s+\epsilon(s)}(K^{(m)})])\no\\
  & \quad \quad \quad \quad \quad\quad \quad \quad \quad \quad\quad \quad \quad \quad \quad\quad \quad \quad \quad \quad\quad \quad \quad -f(s, Y_s, Z_s, Y_{s+\mu(s)},G_{s}(K) ,\mE^{\sF_{s}}[G_{s+\epsilon(s)}(K)])|^{2}\dif s\bigg]\bigg..
  \end{align}
Since $\lim_{n\rightarrow \infty}\mE\bigg[\bigg.\int^{T}_{0}|Z^{(n)}_{s}-Z_{s}|^{2}\dif s\bigg]\bigg. =0,\lim_{n\rightarrow \infty}\mE\bigg[\bigg.\int^{T}_{0}|K^{(n)}_{s}-K_{s}|^{2}\dif s\bigg]\bigg. =0,$  there exists two process $Z'$  in $L^{2}_{\mF}(0,T; \mR^{d}),$ $K'$  in $L^{2}_{\mF}(0,T; \mR),$ and a subsequence if necessary such that $|Z^{(n)}_{t}|\leq Z'_{t}$ and $Z^{(n)}_{t}\rightarrow Z_{t}, |K^{(n)}_{t}|\leq K'_{t}$ and $K^{(n)}_{t}\downarrow K_{t},\dif t\times \dif P-a.e.$
By $(\mathrm{\romannumeral7})$ of Lemma 3.5 and Dominated convergence theorem, we have
$$\lim_{n\rightarrow \infty}\mE\bigg[\bigg.\int^{T}_{0}|f_{n}(s, Y^{(n)}_s, Z^{(n)}_s, Y^{(n)}_{s+\mu(s)},G_{s}(K^{(n)}), \mE^{\sF_{s}}[G_{s+\epsilon(s)}(K^{(n)})])-f(s, Y_s, Z_s, Y_{s+\mu(s)},G_{s}(K), \mE^{\sF_{s}}[G_{s+\epsilon(s)}(K)])|^{2}\dif s\bigg]\bigg.=0.$$
Thus, $(Y, Z, K)\in S^{2}_{\mF}(0, T+\delta; \mR)\times L^{2}_{\mF}(0, T+\delta; \mR^{d}) \times H^{2}_{\mF}(0, T+\delta; \mR)$ satisfying Eq.(6).
Assume that
$(\tilde{Y}, \tilde{Z},\tilde{K})$ is another solution of Eq.(6). Since $f_{n}(s,y, z, \theta_{r},m,\bar{m})\leq f(s, y, z, \theta_{r},m,\bar{m})$ for any $t\in [0, T],r\in [t, T+K],y,m,\bar{m}\in \mR, z\in \mR^{d},\theta\in L^{2}_{\sF}(t, T+K; \mR),$  by the comparison theorems in Theorem 3.6, we have $ Y^{(n)}_{t}\leq \tilde{Y}_{t}, \tilde{K}^{(n)}_{t}\geq \tilde{K}_{t}  ,a.e., a.s.,$  since $Y^{(n)}_{t}\uparrow Y_{t}, n\rightarrow \infty,$ we have  $ Y_{t}\leq \tilde{Y}_{t},a.e., a.s.$
\end{proof}

\section*{Funding}
This work are  supported by  the Special Fund of Basic Scientific Research of Central
Colleges (GZQ16022) and NSF of China (No. 11626236).

\section*{Availability of data and materials}
\begin{center}
Not applicable.
\end{center}

\section*{Competing interests}
\begin{center}
The author declare they have no competing interests.
\end{center}

\section*{Authors' contributions}
All authors conceived of the study and participated in its design and coordination. All authors read and approved the final manuscript.

\end {document}